\documentclass[11pt]{article}
\usepackage[letterpaper,margin=1in]{geometry}
\usepackage{times}
\usepackage{amsmath,amssymb,amsthm}
\usepackage{graphicx}
\usepackage{authblk}
\usepackage{titlesec}
\usepackage{enumitem}
\usepackage[colorlinks=true,linkcolor=red,citecolor=red,urlcolor=blue]{hyperref}

\titleformat{\section}{\large\bfseries}{\thesection.}{0.5em}{}
\titleformat{\subsection}{\normalsize\bfseries}{\thesubsection}{0.6em}{}

\newtheorem{proposition}{Proposition}
\newtheorem{lemma}{Lemma}
\newtheorem{theorem}{Theorem}

\theoremstyle{definition}
\newtheorem{remark}{Remark}

\newcommand{\gc}{\gamma_{\mathrm{crit}}}

\title{\textbf{Exact Solutions and a Critical-Damping Threshold for a\\
Nonlinear Transport-Wave Equation}}
\author{Sameh Gana\textsuperscript{a} and Salma Ben Khelifa\textsuperscript{b}}
\renewenvironment{abstract}{\begin{center}\small\textbf{Abstract}\end{center}
\list{}{\leftmargin=0.5in\rightmargin=0.5in}\item\relax\small}{\endlist}
\date{}
\begin{document}
\maketitle
{\renewcommand{\thefootnote}{}
	\footnotetext{\textsuperscript{a}Department of Basic Sciences, 
		Deanship of Preparatory Year and Supporting Studies, 
		Imam Abdulrahman Bin Faisal University, P.O. Box 1982, Dammam, 34212, Saudi Arabia. Corresponding author: sbgana@iau.edu.sa.
		}
	\footnotetext{\textsuperscript{b}International French School Victor Hugo, Djerba, Tunisia.}}

\begin{abstract}
Linear stability theory does not determine whether a small disturbance in a physical
system decays or grows without bound when the disturbance has finite size and the damping
varies in time. We address this question for the nonlinear damped transport-wave equation
$u_{tt}+cu_{xt}+\gamma u_t=\lambda u(u_t+cu_x)e^{-\gamma t}$, a constant-coefficient case
of a model introduced in previous work. Combining the method of characteristics with a
Riccati reduction, we reduce the equation exactly to a first-order ordinary differential
equation at each point and solve it explicitly. For a natural class of initial data, the solution has a critical damping rate $\gc$, given in closed form in terms of the initial data: for $\gamma<\gc$ the solution
blows up in finite time, and otherwise no finite-time singularity occurs. A similar result holds for time-dependent damping and for power-law feedback. Applications to shock formation in the damped Burgers
equation, gelation with decaying reactivity and frictional heating under a moving band
are presented.

\medskip
\noindent\textbf{Keywords:} Nonlinear partial differential equations, Method of
characteristics, Riccati equation, Thermal runaway, Critical damping,
Finite-time blow-up
\end{abstract}

\section{Introduction}
When a brake pad briefly contacts a moving surface, frictional heating raises the
temperature, and the surface then cools as the heat dissipates. If the friction rises
with temperature, however, a small temperature rise can increase the heating and lead
to thermal runaway. Linear stability analysis examines small perturbations of a steady
state and determines whether they grow or decay. This describes stability against very
small perturbations, but it does not determine what happens to a disturbance of finite
size, particularly when the driving changes in time, for example as surfaces heat up or
wear, or as a load relaxes. For the class of equations introduced in~\cite{MGA24}, the
behaviour of finite disturbances has been studied through particular numerical
solutions.

In this paper, we study the nonlinear damped transport-wave equation
$$u_{tt}+c\,u_{xt}+\gamma u_t=\lambda u\bigl(u_t+cu_x\bigr)e^{-\gamma t},$$ in which a
quantity $u(x,t)$ is transported at speed $c$, damped at rate $\gamma$, and subject to a
nonlinear feedback whose strength depends on $u$ itself. It is a constant-coefficient
case of a class introduced in our previous work~\cite{MGA24}, building on the
reduction technique in~\cite{HS09} and on our earlier partial
differential equation classes~\cite{MGA23}. That paper is a solution-methods
contribution: it reduces the class to first-order equations, identifies the case
$G(u)=u^n$ with $n=2$ as an Abel equation with reference to the classical
handbooks~\cite{Zw98,Mu60,Ka77,PZ11}, and observes in plotted examples that some
solutions become singular. 
The conditions under which these singularities occur, and the physical meaning of the
equation, remain to be determined. Our main result is a closed-form criterion for this
equation. The method combines characteristics with a Riccati reduction. Along
characteristics, the equation reduces exactly to a first-order differential equation. For a natural class of initial data, defined in Section~\ref{sec:riccati}, the equation
at each fixed point becomes, in a rescaled time, an autonomous Riccati equation. Its solution
blows up in finite time if and only if the damping rate $\gamma$ is below a critical value
$\gc$, which we give in closed form in terms of the initial data (Theorem~\ref{thm:threshold}). For other initial data, the set of points at which blow-up
occurs depends on the initial profile and the transport speed, and extends farther
downstream than upstream (Section~\ref{sec:surfaces}). 
\paragraph{Contributions.} The time change that turns a damped scalar equation into an
autonomous one is classical, and so is the fact that the resulting blow-up is decided by
comparing two times. The contributions of this paper are the following.
\begin{enumerate}[label=(\roman*),leftmargin=2.2em]
	\item An exact reduction of a second-order nonlinear transport-wave equation to a
	first-order equation at each fixed point (Proposition~\ref{prop:reduction}), with explicit
	hypotheses on the transport speed and the data.
	\item A sharp critical damping rate $\gc$ and blow-up time $t^*$ in closed form
	(Theorem~\ref{thm:threshold}), extended to time-dependent damping
	(Theorem~\ref{thm:timedamping}) and to power-law feedback, including the Abel case $n=2$
	(Propositions~\ref{prop:power} and~\ref{prop:abel}).
	\item Explicit two-sided bounds for data with variable $K$ (Proposition~\ref{prop:bounds}),
	and numerical evidence that transport shifts the blow-up region downstream.
	\item A physical derivation of the full equation from classical models of frictional
	heating, and two exact correspondences, with shock formation and with gelation, in which
	the threshold has a direct meaning.
\end{enumerate}

The reduced equation describes exactly two classical problems, shock formation in the
damped Burgers equation and gelation with decaying reactivity, to which the threshold
applies unchanged. The full equation arises from a physical model:
Section~\ref{sec:friction} combines classical models of frictional heating at a sliding
contact~\cite{Ja42,No15}, a Maxwell model of band tension~\cite{PC21} and the adiabatic
approximation of thermal-explosion theory~\cite{FK69}. In this model, $u$ is a
temperature, $\gamma$ is the tension-relaxation rate, and blow-up signals the onset of
thermal runaway. The threshold compares the relaxation time of the band with the
adiabatic induction time.

The closest results come from three fields. For semilinear damped wave equations
$u_{tt}-\Delta u+b(t)u_t=|u|^p$, blow-up is governed by critical exponents in $p$ and by
lifespan estimates~\cite{TY01,Wa14,ISW19,MR23}. For scale-invariant damping $b(t)=\frac{\mu}{1+t}$, small-data solutions exist globally
for $p>1+\frac{2}{n}$ when $\mu$ is large enough~\cite{DA15}. Our threshold is of a different kind: the nonlinearity is of transport
type, and $\gc$ depends on the size of the data, not only on an exponent. For the same
type of damping, $b(t)=\frac{\gamma}{1+t}$, Theorem~\ref{thm:timedamping} gives blow-up if
and only if $\gamma<1+\frac{1}{\tau^*}$, a condition that involves the initial amplitude
through $\tau^*$. For reaction--diffusion equations $u_t=\Delta u+h(t)u^p$, a
time-dependent coefficient also shifts the threshold between blow-up and global
existence, and the threshold is again a critical exponent~\cite{Me90}. In thermal
explosion theory, reactant consumption weakens the source over time and changes the
criteria for explosion~\cite{BGW77}. Our friction model gives such a criterion in closed
form, $\tau_r>\tau^*$, for a transport equation. In~\cite{MGA24}, the reduction was
obtained and some solutions were seen to become singular. Here, we determine exactly when this happens.

\paragraph{Organization of the paper.}
We begin in Section~\ref{sec:model} with the general nonlinear damped transport-wave
model and its exact reduction along characteristics. Section~\ref{sec:riccati} treats the constant-coefficient case, in which the reduced equation is a
Riccati equation. We solve it and obtain the critical damping rate $\gc$ in closed form.
We then extend the result to time-dependent damping and to power-law feedback, including
the Abel case $n=2$, and illustrate it numerically, first at a single point and then for
a localized initial disturbance, for which we also give explicit bounds. Section~\ref{sec:application}
presents three applications of these results: shock formation under linear damping,
gelation with decaying reactivity, and frictional heating under a moving band. The
interpretation and limitations of the results are discussed in
Section~\ref{sec:discussion}, and Section~\ref{sec:conclusion} concludes the paper.

\section{The general nonlinear damped transport-wave model}\label{sec:model}

We now consider the general second-order equation, a special case of the class we
introduced in our previous work~\cite{MGA24}:
\begin{equation}
u_{tt}+a(x,t)u_{xt}+b(t)u_t=G(u)\bigl(u_t+a(x,t)u_x\bigr)e^{-\int_0^t b(s)\,ds},
\label{eq:general}
\end{equation}
where $a(x,t)$ is a transport speed, $b(t)$ a 
time-varying damping rate, and $G(u)$ is an arbitrary nonlinear feedback function. The reduction below reproduces our earlier argument in~\cite{MGA24}. The threshold
results of Section~\ref{sec:riccati} and their applications in
Section~\ref{sec:application} are new to the present paper.

\subsection{Reduction along characteristics}

Assume that:
\begin{enumerate}[label=$ (H_{\arabic*})$,leftmargin=2.6em]
\item $a\in C^1(\Omega_T)$, and
$|a(x,t)|\le M(t)\,(1+|x|)$ for a continuous function $M$.
\item $b\in C([0,T))$, $G\in C^1(\mathbb{R})$ and $g'=G$.
\item  $u(\cdot,0)=u_0\in C^2(\mathbb{R})$ and
$u_t(\cdot,0)=u_1\in C^1(\mathbb{R})$,
\end{enumerate}
where $\Omega_T =\mathbb{R}\times[0,T)$ with $0<T\le\infty$.
\paragraph{Characteristics and the foot map $X_0$.}
Since $a\in C^1(\Omega_T)$, it is locally Lipschitz
with respect to its spatial variable. Hence, by the
Picard--Lindel\"of theorem, for every $x_0\in\mathbb{R}$,
the characteristic equation
\begin{equation}
	\begin{cases}
		\dot{x}(t)=a(x(t),t),\\
		x(0)=x_0
	\end{cases}
	\label{eq:characteristic-ode}
\end{equation}
has a unique solution near $t=0$, which we denote by $x(t;x_0)$. By (H1),
\[
1+|x(t;x_0)|
\leq 1+|x_0|
+\int_0^t M(s)(1+|x(s;x_0)|)\,ds.
\]
By Gr\"onwall's inequality,
\begin{equation}
	1+|x(t;x_0)|
	\leq (1+|x_0|)
	\exp\left(\int_0^t M(s)\,ds\right).
	\label{eq:characteristic-bound}
\end{equation}
Since this bound is finite for $t<T$, the characteristic exists on all of $[0,T)$. The
same argument backward in time shows that the characteristic through any point of
$\Omega_T$ reaches $t=0$. For each $t\in[0,T)$, we define the characteristic flow
\[
\Phi_t:\mathbb{R}\longrightarrow\mathbb{R},
\quad
\Phi_t(x_0)=x(t;x_0).
\]
By uniqueness, $\Phi_t$ is injective, and by the backward argument it is onto, so $\Phi_t$
is a bijection of $\mathbb{R}$. Differentiating~\eqref{eq:characteristic-ode} with respect
to $x_0$ gives
\begin{equation}
	\frac{\partial x(t;x_0)}{\partial x_0}
	=\exp\left(\int_0^t a_x\bigl(x(s;x_0),s\bigr)\,ds\right)>0.
	\label{eq:flow-derivative}
\end{equation}
For each $(x,t)\in\Omega_T$, let $X_0(x,t)=\Phi_t^{-1}(x)$ be the starting point of the characteristic
through $(x,t)$. Under (H1), $X_0\in C^1(\Omega_T)$ and $X_0$ is constant along
characteristics, $\partial_tX_0+a\,\partial_xX_0=0$. For $a=c$, used in the rest of the paper,
$X_0(x,t)=x-ct$.
By (H2) and (H3), $C=u_1-g(u_0) \in C^1(\mathbb{R})$. On the characteristic starting at $x_0$, we have
$C\bigl(X_0(x,t)\bigr)=C(x_0)$, and
\begin{equation}
	\bigl(\partial_t+a\,\partial_x\bigr)\,C\bigl(X_0(x,t)\bigr)=0 .
	\label{eq:CX0}
\end{equation}

\begin{proposition}[Reduction along characteristics]\label{prop:reduction}
Assume (H1), (H2) and (H3), and let $u\in C^2(\Omega_T)$ be a solution of~\eqref{eq:general}
with $u(\cdot,0)=u_0$ and $u_t(\cdot,0)=u_1$. Then, along each characteristic
$x(t)=x(t;x_0)$, the function $P(t)=u_t(x(t),t)$ satisfies
\begin{equation}
e^{\int_0^t b(s)\,ds}\,P(t)=g\bigl(U(t)\bigr)+C(x_0),\qquad U(t)=u(x(t),t),
\label{eq:prop}
\end{equation}
for a constant $C(x_0)=u_1(x_0)-g\bigl(u_0(x_0)\bigr)$ which depends only on the characteristic. In other words, ~\eqref{eq:general} reduces to the
first-order equation
\begin{equation}
u_t(x,t)=e^{-\int_0^t b(s)\,ds}\Bigl[g\bigl(u(x,t)\bigr)+C\bigl(X_0(x,t)\bigr)\Bigr],
\label{eq:reduced}
\end{equation}
which, at each fixed point $x$, is an ordinary
differential equation in $t$. Conversely, let $C\in C^1(\mathbb{R})$ and let $u\in C^2(\Omega_T)$ satisfy~\eqref{eq:reduced}
on $\Omega_T$. Then $u$ solves~\eqref{eq:general} on $\Omega_T$ and
$C=u_t(\cdot,0)-g\bigl(u(\cdot,0)\bigr)$.
\end{proposition}

\begin{proof}
Along the characteristic, $U'(t)=u_t+a(x,t)u_x$ and
$P'(t)=u_{tt}+a(x,t)u_{xt}$. Equation~\eqref{eq:general}
gives
\begin{equation}
P'(t)+b(t)P(t)=G\bigl(U(t)\bigr)\,U'(t)\,e^{-\int_0^t b(s)\,ds}.
\end{equation}
Multiplying both sides by $\mu(t)=e^{\int_0^t b(s)\,ds}$ gives
\begin{equation}
\frac{d}{dt}\bigl[\mu(t)P(t)\bigr]=G\bigl(U(t)\bigr)\,U'(t)=\frac{d}{dt}\bigl[g\bigl(U(t)\bigr)\bigr],
\end{equation}
Integrating from $0$ to $t$ gives $\mu(t)P(t)=g(U(t))+C(x_0)$, which is ~\eqref{eq:prop}. Since $P(t)=u_t(x(t),t)$ is the value of $u_t$ at the point $(x(t),t)$, and the
characteristic through $(x,t)$ starts at $x_0=X_0(x,t)$, this proves~\eqref{eq:reduced}.\\
We now prove the second statement. Assume that $u\in C^2(\Omega_T)$ satisfies~\eqref{eq:reduced}. Then
$\mu(t)\,u_t-g(u)=C\bigl(X_0(x,t)\bigr)$. Applying $\partial_t+a\,\partial_x$ and
using~\eqref{eq:CX0}, $\mu'=b\mu$ and $g'=G$ gives 
\[
\mu(t)\bigl(u_{tt}+a\,u_{xt}+b\,u_t\bigr)-G(u)\bigl(u_t+a\,u_x\bigr)=0 .
\]
Dividing by $\mu(t)>0$ gives~\eqref{eq:general}. At $t=0$, $\mu(0)=1$ and $X_0(x,0)=x$,
so~\eqref{eq:reduced} gives \[u_t(x,0)=g\bigl(u(x,0)\bigr)+C(x).\]
\end{proof}
\begin{remark}
Equation~\eqref{eq:reduced} involves $u_t$ at a fixed point, not the derivative
$u_t+a\,u_x$ along a characteristic. At a fixed point $x$, different characteristics
pass through $x$ over time, so the forcing $C(X_0(x,t))$ changes with $t$ unless $C$ is
constant.
\end{remark}
In Proposition~\ref{prop:reduction}, we reduce the second-order equation~\eqref{eq:general}
exactly to a first-order ordinary differential equation at each fixed point. For constant
coefficients, $G(u)=\lambda u$ and constant $C$, this equation is a Riccati equation,
which we solve explicitly in Section~\ref{sec:riccati}.

\section{The constant-coefficient case}\label{sec:riccati}

Let $a(x,t)= c$, $b(t)=\gamma$ and $G(u)=\lambda u$.
Equation~\eqref{eq:general} takes the form
\begin{equation}
u_{tt}+c\,u_{xt}+\gamma u_t=\lambda u\bigl(u_t+cu_x\bigr)e^{-\gamma t}.
\label{eq:model}
\end{equation}
Then ~\eqref{eq:reduced} reduces to 
\begin{equation}
u_t(x,t)=e^{-\gamma t}\Bigl[\tfrac{\lambda}{2}u(x,t)^2+K(x-ct)\Bigr],\qquad
K(\xi)=u_t(\xi,0)-\tfrac{\lambda}{2}u(\xi,0)^2.\
\label{eq:reducedmodel}
\end{equation}
In Sections~\ref{sec:removing}--\ref{sec:example} we take initial data for which $K$ is
the same constant on every characteristic, i.e.\
$u_t(x,0)=K+\frac{\lambda}{2}u(x,0)^2$ for all $x$. Then~\eqref{eq:reducedmodel} holds at every point $x$, with only the initial value
$u(x,0)$ depending on $x$. The case where $K$ varies in
space is treated in Section~\ref{sec:surfaces}.

\subsection{Removing time explicitly}\label{sec:removing}

Let $ t\in[0,\infty)$ and 
\begin{equation}
\tau(t)=\frac{1-e^{-\gamma t}}{\gamma},\qquad \tau\in\Bigl[0,\frac{1}{\gamma}\Bigr)
,
\label{eq:tau}
\end{equation}
 If $V(\tau)=u(x,t(\tau))$ at a fixed point $x$,
then we obtain a Riccati equation
\begin{equation}
\frac{dV}{d\tau}=\frac{\lambda}{2}V^2+K.
\label{eq:riccati}
\end{equation}
If $A=u(x,0)$, $B=u_t(x,0)$, then 
$K=B-\frac{\lambda}{2}A^2$.

\begin{lemma}\label{lem:riccati}
Let $\lambda>0$ and $V$ be the solution of~\eqref{eq:riccati} with $V(0)=A$.
\begin{enumerate}[label=(\roman*),leftmargin=2.2em]
\item For \text{$K=k^2>0$,}
$V(\tau)=k\sqrt{\frac{2}{\lambda}}\,\tan\Bigl[k\sqrt{\frac{\lambda}{2}}\,(\tau+\tau_0)\Bigr]$,
with $\tau_0$ fixed by $V(0)=A$.
\item For \text{$K=-k^2<0$ and $|A|<k\sqrt{\frac{2}{\lambda}}$, }
$V(\tau)=-k\sqrt{\frac{2}{\lambda}}\,\tanh\Bigl[k\sqrt{\frac{\lambda}{2}}\,(\tau-\tau_0)\Bigr]$, with $\tau_0$ fixed by $V(0)=A$, as in (i).
\item For \text{$K=0$}, $V(\tau)=\dfrac{A}{1-\frac{\lambda}{2}A\tau}$.
\end{enumerate}
\end{lemma}

\begin{proof}
	Separating variables in~\eqref{eq:riccati} gives
	\[
	\tau=\int_A^{V}\frac{ds}{\frac{\lambda}{2}s^2+K},
	\] where $V=V(\tau)$.\\
	For $K=k^2>0$,
	\[
	\tau=\frac{1}{k\sqrt{\frac{\lambda}{2}}}\Bigl[\arctan\frac{V}{k\sqrt{\frac{2}{\lambda}}}-\arctan\frac{A}{k\sqrt{\frac{2}{\lambda}}}\Bigr].
	\]
	For $K=-k^2<0$ and $|A|<k\sqrt{\frac{2}{\lambda}}$,
	\[
	\tau=\frac{1}{k\sqrt{\frac{\lambda}{2}}}\Bigl[\operatorname{artanh}\frac{A}{k\sqrt{\frac{2}{\lambda}}}-\operatorname{artanh}\frac{V}{k\sqrt{\frac{2}{\lambda}}}\Bigr].
	\]
	For $K=0$,
	\[
	\tau=\frac{2}{\lambda}\Bigl(\frac{1}{A}-\frac{1}{V}\Bigr).
	\]
	Solving each relation for $V$ completes the proof.
\end{proof}

\subsection{The blow-up condition}

The tan branch diverges at a finite value $\tau^*$ given by
\begin{equation}
	\tau^*=\frac{1}{k\sqrt{\frac{\lambda}{2}}}\left(\frac{\pi}{2}-\arctan\frac{A}{k\sqrt{\frac{2}{\lambda}}}\right).
	\label{eq:taustar}
\end{equation}

\begin{theorem}[Critical damping rate]\label{thm:threshold}
	Let $a=c$, $b=\gamma>0$ and $G(u)=\lambda u$ with $\lambda>0$. Assume that $K$ is constant, with $K=k^2>0$. For each point $x$, let $\tau^*$ be given
	by~\eqref{eq:taustar} with $A=u(x,0)$ and define the critical damping rate
	$\gc=1/\tau^*$.
\begin{enumerate}[label=(\alph*),leftmargin=2.2em]
	\item If $\gamma<\gc$, there is finite-time blow-up at
	\begin{equation}
		t^*=-\frac{1}{\gamma}\ln\bigl(1-\gamma\tau^*\bigr).
		\label{eq:threshold}
	\end{equation}
	\item If $\gamma=\gc$, there is no finite-time singularity, but $u\to\infty$ as
	$t\to\infty$.
	\item If $\gamma>\gc$, there is no finite-time singularity, and $u$ approaches a finite
	limit as $t\to\infty$.
\end{enumerate}
\end{theorem}

\begin{proof}
	Using Lemma~\ref{lem:riccati} (i), $V(\tau)$ is finite for $\tau<\tau^*$ and
	$V(\tau)\to\infty$ as $\tau\to\tau^*$ and $$u(x,t)=V(\tau(t)),$$ where $\tau(t)=\frac{1-e^{-\gamma t}}{\gamma}$ increases and
	tends to $\frac{1}{\gamma}$ without reaching it.
	\begin{enumerate}[label=(\alph*),leftmargin=2.2em]
		\item If $\gamma<\gc$, then $\tau^*<\frac{1}{\gamma}$, and $\tau(t)$ reaches $\tau^*$ at a
		finite time. Solving $\tau(t^*)=\tau^*$ gives~\eqref{eq:threshold}.
		\item If $\gamma=\gc$, then $\tau(t)\to\tau^*$ only as $t\to\infty$, so $u$ becomes
		infinite, but not in finite time.
		\item If $\gamma>\gc$, then $\tau(t)\to\frac{1}{\gamma}<\tau^*$, so
		$u\to V\bigl(\frac{1}{\gamma}\bigr)<\infty$.
	\end{enumerate}
\end{proof}
\begin{remark}
	Let $K=-k^2<0$. If $A\le k\sqrt{\frac{2}{\lambda}}$, then $V$ stays bounded. If
	$A>k\sqrt{\frac{2}{\lambda}}$, then $V$ blows up at
	\[
	\tau^*=\frac{1}{k\sqrt{\frac{\lambda}{2}}}\,\operatorname{arcoth}\frac{A}{k\sqrt{\frac{2}{\lambda}}},
	\]
	and the conclusions of Theorem~\ref{thm:threshold} hold with $\gc=1/\tau^*$.
\end{remark}
\begin{remark}
	Theorem~\ref{thm:threshold} holds at each point $x$ separately. Since $\tau^*$ decreases
	as $A$ increases, $\gc$ increases with $u(x,0)$. Since
	$u_t(x,0)=K+\frac{\lambda}{2}u(x,0)^2$, the condition $u(\cdot,0)\in C^2(\mathbb{R})$ gives
	$u_t(\cdot,0)\in C^2(\mathbb{R})$, and $u(x,t)=V\bigl(\tau(t)\bigr)$ is a $C^2$ solution on
	the set $\{(x,t):t<t^*(x)\}$, where $t^*(x)=\infty$ at points that do not blow up, since
	each point evolves independently. The classical solution therefore exists on
	$\mathbb{R}\times[0,T)$ with $T=\inf_x t^*(x)$, and no point blows up if and only if
	$\gamma\ge\gc\bigl(\sup_x u(x,0)\bigr)$. If $u(\cdot,0)$ is unbounded above, then $t^*(x)$
	takes arbitrarily small values, so $T=0$ for every $\gamma$. For $K<0$ and $K=0$,
	Remark~2 and Proposition~\ref{prop:power} give the same conclusions.
\end{remark}

\subsection{Time-dependent damping}\label{sec:timedamping}

Theorem~\ref{thm:threshold} extends to a time-dependent damping rate $b(t)$. Let $a=c$,
$G(u)=\lambda u$ with $\lambda>0$, $b\in C([0,\infty))$, and let $K=k^2>0$ be the same on
every characteristic. By Proposition~\ref{prop:reduction}, at each fixed point $x$,
\[
u_t=e^{-\int_0^t b(s)\,ds}\Bigl[\frac{\lambda}{2}u^2+K\Bigr].
\]
Set
\[
\tau(t)=\int_0^t e^{-\int_0^s b(r)\,dr}\,ds,\qquad
\tau_\infty=\lim_{t\to\infty}\tau(t)\in(0,\infty].
\]
Then $\tau$ increases strictly from $0$ to $\tau_\infty$, and $V(\tau)=u(x,t)$ satisfies
the same Riccati equation~\eqref{eq:riccati}.

\begin{theorem}\label{thm:timedamping}
	Under these assumptions, let $\tau^*$ be given by~\eqref{eq:taustar} with $A=u(x,0)$.
	\begin{enumerate}[label=(\alph*),leftmargin=2.2em]
		\item If $\tau^*<\tau_\infty$, there is finite-time blow-up at the time $t^*$ defined
		by $\tau(t^*)=\tau^*$.
		\item If $\tau^*=\tau_\infty<\infty$, there is no finite-time singularity, but
		$u\to\infty$ as $t\to\infty$.
		\item If $\tau^*>\tau_\infty$, there is no finite-time singularity, and $u$ approaches
		the finite limit $V(\tau_\infty)$ as $t\to\infty$.
	\end{enumerate}
\end{theorem}

\begin{proof}
	Since $u(x,t)=V(\tau(t))$, where $\tau$ is continuous and increasing, the proof is that of
	Theorem~\ref{thm:threshold} with $\frac{1}{\gamma}$ replaced by $\tau_\infty$.
\end{proof}

For constant damping, $\tau_\infty=\frac{1}{\gamma}$ and Theorem~\ref{thm:timedamping}
reduces to Theorem~\ref{thm:threshold}. If $b$ is integrable on $[0,\infty)$, for example
if the damping decays exponentially, then $\tau_\infty=\infty$ and the solution always
blows up. For the decaying rate $b(t)=\frac{\gamma}{1+t}$, we have
$e^{-\int_0^t b}=(1+t)^{-\gamma}$, so $\tau_\infty=\frac{1}{\gamma-1}$ for $\gamma>1$ and
$\tau_\infty=\infty$ for $\gamma\le1$. The solution therefore blows up in finite time if
and only if $\gamma<1+\frac{1}{\tau^*}$.

\subsection{Generalization: power-law feedback $G(u)=\lambda u^n$}

We now replace $G(u)=\lambda u$ of~\eqref{eq:model} by the power law
$G(u)=\lambda u^n$, which reduces to it for $n=1$. In our previous work~\cite{MGA24}, we showed that for $n=2$ the general
equation~\eqref{eq:general} leads to an Abel equation, whose solution is discussed in
the classical handbooks~\cite{Zw98,Mu60,Ka77,PZ11}. This
subsection derives the threshold for the power-law family, and Section~\ref{sec:abel}
treats the case $n=2$ explicitly.

With $G(u)=\lambda u^n$, the antiderivative in Proposition~\ref{prop:reduction} is
$g(u)=\frac{\lambda}{n+1}u^{n+1}$, so the reduced equation~\eqref{eq:riccati}
generalizes to
\begin{equation}
\frac{dV}{d\tau}=\frac{\lambda}{n+1}V^{n+1}+K.
\label{eq:abel}
\end{equation}
For $n=1$, ~\eqref{eq:abel} is the Riccati equation of Lemma~\ref{lem:riccati}.
For $n=2$, it is an Abel equation of the first kind (Section~\ref{sec:abel}). For every
$n$, the case $K=0$, which corresponds to $B=\frac{\lambda}{n+1}A^{n+1}$, is solved explicitly by
separation of variables.

\begin{proposition}\label{prop:power}
Let $G(u)=\lambda u^n$ with $\lambda>0$, $n\ge1$ an integer, $A>0$, and $K=0$. Then
\begin{equation}
V(\tau)=\Bigl[A^{-n}-\frac{n\lambda}{n+1}\,\tau\Bigr]^{-1/n},
\label{eq:Vpower}
\end{equation}
which blows up at the rescaled time $$\tau^*=\frac{n+1}{n\lambda A^n},$$ and the
conclusions of Theorem~\ref{thm:threshold} hold with
\begin{equation}
\gc(n)=\frac{n\lambda A^n}{n+1}.
\label{eq:gcn}
\end{equation}
\end{proposition}

\begin{proof}
	For $K=0$, equation~\eqref{eq:abel} gives $V'=\frac{\lambda}{n+1}V^{n+1}$ with $V(0)=A>0$.
	Using separation of variables,
	\[
	\int_A^{V(\tau)}s^{-(n+1)}\,ds=\frac{\lambda}{n+1}\,\tau,
	\qquad\text{that is,}\qquad
	\frac{A^{-n}-V(\tau)^{-n}}{n}=\frac{\lambda}{n+1}\,\tau .
	\]
	Then,
	\[
	V(\tau)^{-n}=A^{-n}-\frac{n\lambda}{n+1}\,\tau ,
	\]
	which leads to~\eqref{eq:Vpower}. The right-hand side decreases from $A^{-n}>0$ and vanishes
	at $\tau^*=\frac{n+1}{n\lambda A^n}$. As $\tau\to\tau^*$, $V(\tau)^{-n}\to0$, so
	$V(\tau)\to\infty$. Then, the conclusions in the original time $t$ follow from the proof of
	Theorem~\ref{thm:threshold}, with $\gc(n)=1/\tau^*$.
\end{proof}
\begin{remark}
	Let $K=0$. If $n$ is even and $A<0$, then $V$ decreases and tends to $-\infty$ at
	$\tau^*=\frac{n+1}{n\lambda|A|^n}$, and the same conclusions hold with $u\to-\infty$. If
	$A=0$, or if $n$ is odd and $A<0$, then $V$ stays bounded.
\end{remark}

\begin{figure}[t]
\centering
\includegraphics[width=0.72\textwidth]{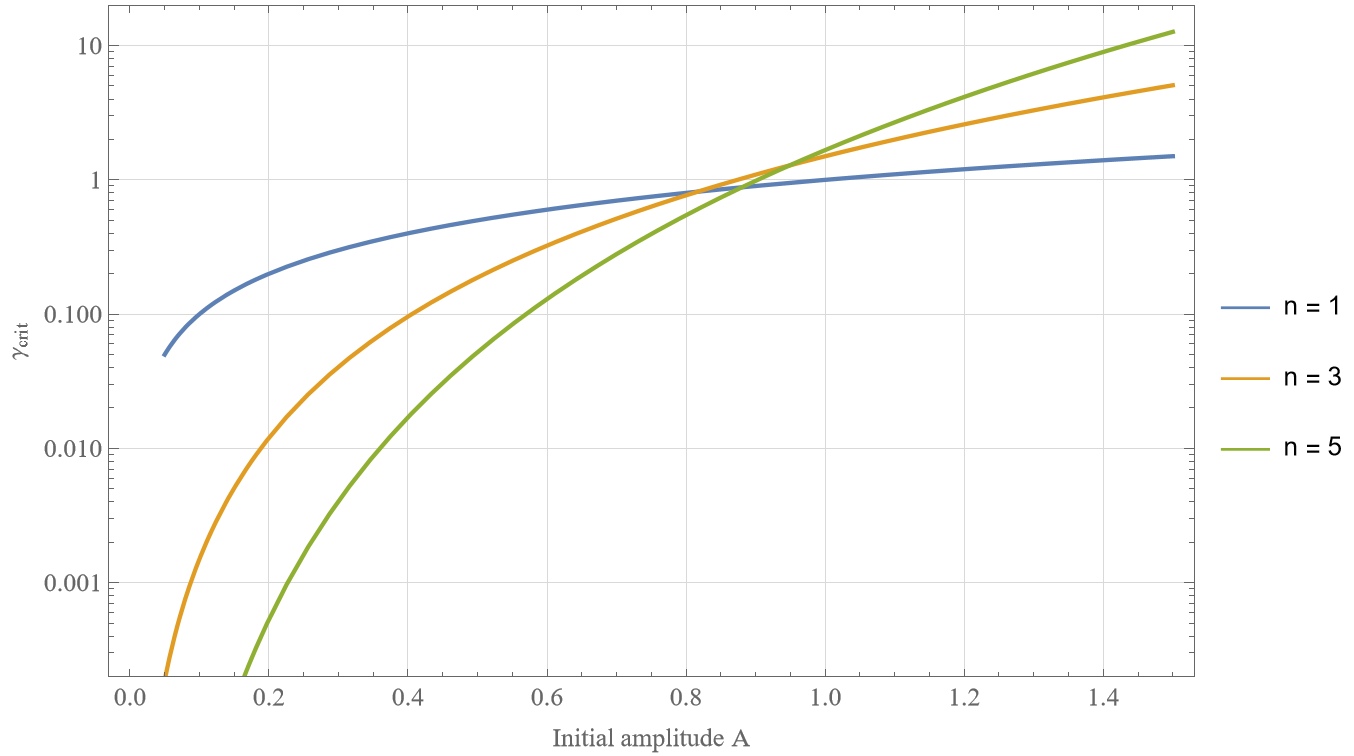}
\caption{Critical damping rate $\gc(n)$ from~\eqref{eq:gcn} as a function of the initial
	amplitude $A$, for $n=1$, $3$ and $5$, with $\lambda=2$.}
\label{fig:gcn}
\end{figure}

\begin{remark}
	For the same value of $\lambda$, $A>0$ and $n\ge2$, we have $\gc(n)<\gc(1)$ if and only if
	\[
	A<A_n=\Bigl(\frac{n+1}{2n}\Bigr)^{\frac{1}{n-1}},
	\]
	where $A_n<1$. Thus higher-order feedback requires less damping for small amplitudes and
	more for large ones (Figure~\ref{fig:gcn}).
\end{remark}

\paragraph{Example $n=3$: cubic feedback.}
Cubic nonlinearities arise in many classical models~\cite{NM79,GH83}.
For $G(u)=\lambda u^3$ with $\lambda=2$, $A=0.5$ and $K=0$, equation~\eqref{eq:gcn} gives
\[
\gc(3)=\frac{3\cdot2\cdot0.5^3}{4}=0.1875,
\]
while linear feedback with the same $\lambda$ and $A$ gives $\gc(1)=\frac{\lambda A}{2}=0.5$.
At this amplitude, cubic feedback needs less damping than linear feedback. Above the
crossover amplitude $A_3=\sqrt{\frac{2}{3}}\approx0.816$, cubic feedback needs more
damping than linear feedback, and the difference grows rapidly with $A$. At $A=1.5$,
$\gc(3)\approx5.06$ while $\gc(1)=1.5$ (Figure~\ref{fig:gcn}).

\subsection{Solving the $K\neq0$ Abel case for $n=2$ by classical partial fractions}
\label{sec:abel}

For $n=2$, equation~\eqref{eq:abel} reads
\[
\frac{dV}{d\tau}=\beta V^3+K,\qquad \beta=\frac{\lambda}{3},
\]
which is an Abel equation of the first kind~\cite{MGA24}. \\Writing $\alpha=\bigl(-\frac{K}{\beta}\bigr)^{1/3}$, the real cube root, the cubic factors as
\begin{equation}
	\beta V^3+K=\beta(V-\alpha)\bigl(V^2+\alpha V+\alpha^2\bigr),
\end{equation}
and the partial-fraction identity~\cite{Zw98,Ka77}
\begin{equation}
	\frac{1}{V^3-\alpha^3}=\frac{1}{3\alpha^2}\Bigl[\frac{1}{V-\alpha}-\frac{V+2\alpha}{V^2+\alpha V+\alpha^2}\Bigr]
\end{equation}
reduces the integral to elementary functions.

\begin{proposition}\label{prop:abel}
	Let $G(u)=\lambda u^2$ with $\lambda>0$, $K\neq0$, $\beta=\frac{\lambda}{3}$ and
	$\alpha=\bigl(-\frac{K}{\beta}\bigr)^{1/3}$. If $A>\alpha$, the solution of
	$$\frac{dV}{d\tau}=\beta V^3+K,V(0)=A$$ is increasing and satisfies $\tau=F(V)-F(A)$, where
	\begin{equation}
		F(V)=\frac{1}{3\beta \alpha^2}\Bigl[\ln|V-\alpha|-\frac{1}{2}\ln\bigl(V^2+\alpha V+\alpha^2\bigr)
		-\sqrt{3}\arctan\frac{2V+\alpha}{\alpha\sqrt{3}}\Bigr].
		\label{eq:F}
	\end{equation}
	It blows up at
	\begin{equation}
		\tau^*=\frac{1}{3\beta \alpha^2}\Bigl[-\frac{\sqrt{3}\,\pi}{2}\,\mathrm{sgn}(\alpha)-\ln|A-\alpha|
		+\frac{1}{2}\ln\bigl(A^2+\alpha A+\alpha^2\bigr)+\sqrt{3}\arctan\frac{2A+\alpha}{\alpha\sqrt{3}}\Bigr],
		\label{eq:taustarabel}
	\end{equation}
	and the conclusions of Theorem~\ref{thm:threshold} hold with $\gc=\frac{1}{\tau^*}$.
\end{proposition}

\begin{proof}
	For $V>\alpha$, we have $\beta V^3+K=\beta(V^3-\alpha^3)>0$.
	Using separation of variables and the partial-fraction identity, we have
	\[
	\tau=\int_A^{V(\tau)}\frac{ds}{\beta(s^3-\alpha^3)}=F\bigl(V(\tau)\bigr)-F(A).
	\]
	As $V\to\infty$, the two logarithmic terms of $F$ cancel, since
	$\frac{V-\alpha}{\sqrt{V^2+\alpha V+\alpha^2}}\to1$, and
	$\arctan\frac{2V+\alpha}{\alpha\sqrt{3}}\to\frac{\pi}{2}\,\mathrm{sgn}(\alpha)$. Then
	\[
	F(V)\to-\frac{\sqrt{3}\,\pi}{6\beta \alpha^2}\,\mathrm{sgn}(\alpha),
	\]
	which is finite. The solution blows up at 
	\[
	\tau^*=-\frac{\sqrt{3}\,\pi}{6\beta \alpha^2}\,\mathrm{sgn}(\alpha)-F(A),
	\] which
	gives ~\eqref{eq:taustarabel}. The conclusions in the original time $t$ then follow from the
	proof of Theorem~\ref{thm:threshold}, with $\gc=\frac{1}{\tau^*}$.
\end{proof}
\begin{remark}
	If $A<\alpha$, then $\frac{dV}{d\tau}=\beta(V^3-\alpha^3)<0$, so $V$ decreases and tends to
	$-\infty$ at the finite rescaled time
	\[
	\tau^*_-=\frac{\sqrt{3}\,\pi}{6\beta \alpha^2}\,\mathrm{sgn}(\alpha)-F(A),
	\]
	obtained as in the proof by letting $V\to-\infty$ in~\eqref{eq:F}. The conclusions of
	Theorem~\ref{thm:threshold} then hold with $\gc=1/\tau^*_-$ and $u\to-\infty$. If $A=\alpha$,
	then $V\equiv \alpha$.
\end{remark}
Equation~\eqref{eq:F} gives $\tau$ as a function of $V$, but it cannot be inverted in closed form. The blow-up time $\tau^*$, however, is explicit, so the
threshold $\gc=\frac{1}{\tau^*}$ is given in closed form for $n=2$ and $K\neq0$.

\subsection{Numerical example}\label{sec:example}

\begin{figure}[t]
\centering
\includegraphics[width=0.85\textwidth]{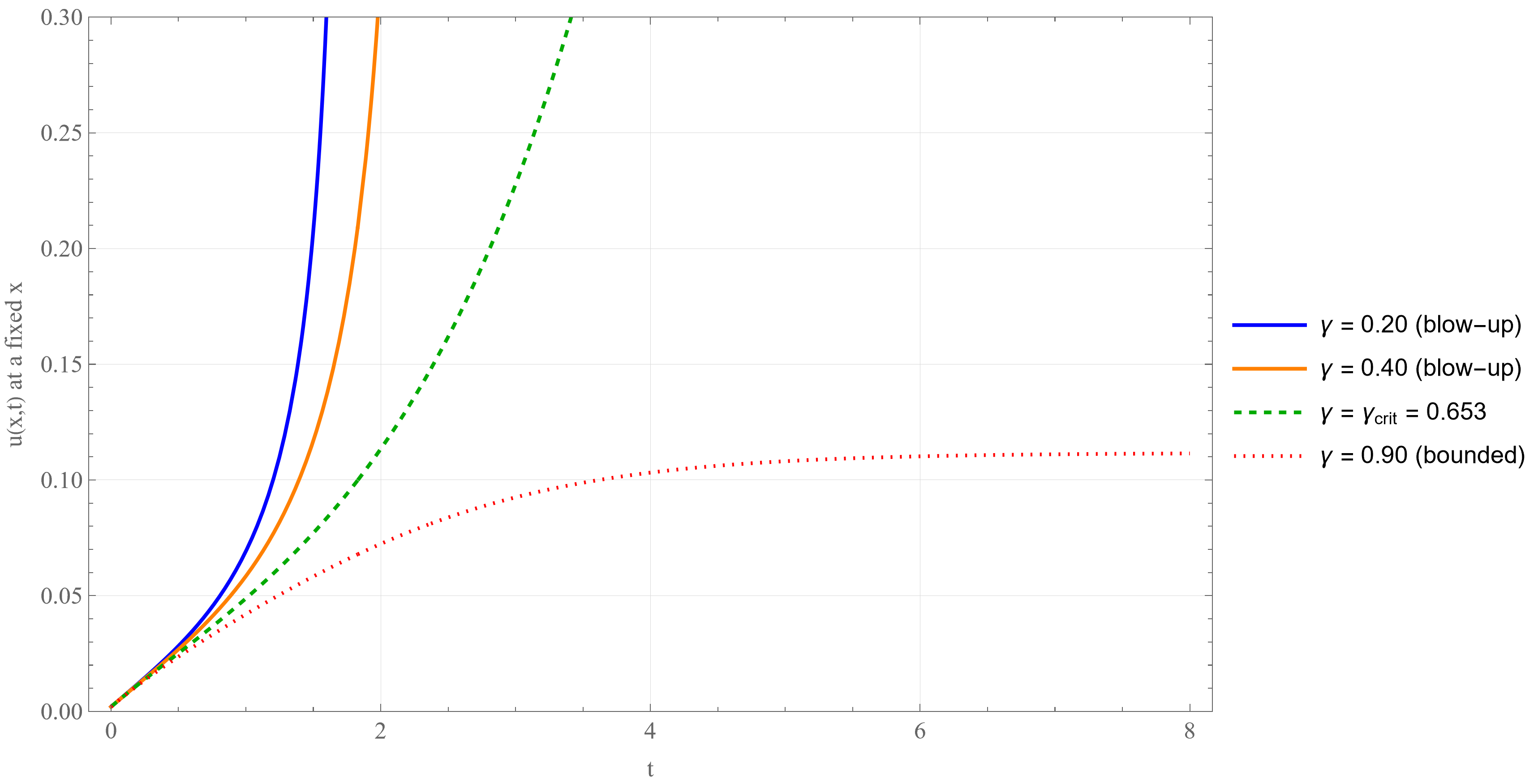}
\caption{Exact solution of equation~\eqref{eq:model} at a fixed point for four
damping rates $\gamma$, with $\lambda=40$, $A=0.002$, $B=0.05$.}
\label{fig:dichotomy}
\end{figure}

Let $\lambda=40$, and take initial data with $A=u(x,0)=0.002$ and
$B=u_t(x,0)=0.05$, with the same value of $K$ on every characteristic. These values are
illustrative and not fitted to a physical system. Then
\[
K=B-\frac{\lambda}{2}A^2=0.0499,\qquad k=\sqrt{K}=0.2234,
\]
and~\eqref{eq:taustar} gives $\tau^*=1.532$. By Theorem~\ref{thm:threshold},
\[
\gc=\frac{1}{\tau^*}=0.6527.
\]
Figure~\ref{fig:dichotomy} shows the exact solution at this point for four damping
rates. For $\gamma=0.2$ and $\gamma=0.4$, below $\gc$, the solution blows up at
$t^*=1.83$ and $t^*=2.37$, given by~\eqref{eq:threshold}. For $\gamma=\gc$, it grows
without bound but never in finite time. For $\gamma=0.9$, above $\gc$, it tends to a
finite limit.

\subsection{A localized initial disturbance}\label{sec:surfaces}
When $K$ is not constant, the forcing $K(x-ct)$ in~\eqref{eq:reducedmodel} varies in space
and is transported at speed $c$. Take $\lambda=40$, $c=1$ and the Gaussian initial data
\begin{equation}
	u(x,0)=0.03\,e^{-x^2/4.5},\qquad u_t(x,0)=0.4\,e^{-x^2/4.5}.
	\label{eq:gauss}
\end{equation}
Then $K(\xi)=u_t(\xi,0)-\frac{\lambda}{2}u(\xi,0)^2$ is not
constant. At a fixed point $x$, the forcing $K(x-ct)$ changes as the pulse moves, and
Theorem~\ref{thm:threshold} does not apply pointwise. Using
Proposition~\ref{prop:reduction}, $u(x,\cdot)$ is the solution
of~\eqref{eq:reducedmodel} at each fixed $x$, which we integrate numerically with
\textsc{Mathematica} (\texttt{NDSolve}). A direct integration of~\eqref{eq:model} by the
method of lines agrees with it to a relative difference of $3.6\times10^{-6}$.

When $K$ varies in space, Theorem~\ref{thm:threshold} still gives two-sided bounds. For
$A>0$ and $K\ge0$, write $\gc(A,K)$ for the critical rate of Theorem~\ref{thm:threshold}
when $K>0$, and $\gc(A,0)=\frac{\lambda A}{2}$ (Lemma~\ref{lem:riccati}\,(iii)).

\begin{proposition}\label{prop:bounds}
Let $a=c$, $b=\gamma>0$, $G(u)=\lambda u$ with $\lambda>0$, assume
$u(\cdot,0),u_t(\cdot,0)\in C^2(\mathbb{R})$, and let $K$ be continuous with
$0\le K_-\le K(\xi)\le K_+$ for all $\xi$.
	\begin{enumerate}[label=(\alph*),leftmargin=2.2em]
		\item If $u(x,0)>0$ and $\gamma<\gc\bigl(u(x,0),K_-\bigr)$, the solution blows up at
		$x$ no later than the blow-up time~\eqref{eq:threshold} computed with
		$A=u(x,0)$ and $K=K_-$.
		\item If $A_{\max}=\sup_x u(x,0)>0$ and $\gamma\ge\gc(A_{\max},K_+)$ with $K_+>0$,
		no point blows up in finite time.
	\end{enumerate}
\end{proposition}

\begin{proof}
	At a fixed $x$, $u$ solves $u_t=e^{-\gamma t}\bigl[\frac{\lambda}{2}u^2+K(x-ct)\bigr]$.
	Let $W_\pm$ solve $W_t=e^{-\gamma t}\bigl[\frac{\lambda}{2}W^2+K_\pm\bigr]$. Since
	$K_-\le K(x-ct)\le K_+$, the comparison theorem for scalar differential equations gives
	$W_-\le u$ when $W_-(0)=u(x,0)$, and $u\le W_+$ when $W_+(0)=A_{\max}$, as long as the
	functions exist. In (a), $W_-$ blows up by Theorem~\ref{thm:threshold}, or by
	Proposition~\ref{prop:power} if $K_-=0$, so $u$ does too. In (b), $u_t\ge0$ because
	$K\ge0$, so $u(x,0)\le u\le W_+$, and $W_+$ has no finite-time singularity by
	Theorem~\ref{thm:threshold}.
\end{proof}

For the data~\eqref{eq:gauss}, $K(\xi)=0.4E-20(0.03E)^2$ with $E=e^{-\xi^2/4.5}\in(0,1]$,
so $0<K(\xi)\le K(0)=0.382$, and we may take $K_-=0$ and $K_+=0.382$. At $x=0$,
$u(0,0)=0.03$ and $\gc(0.03,0)=\frac{\lambda A}{2}=0.6$, so the solution blows up at $x=0$
for every $\gamma<0.6$. With $A_{\max}=0.03$ and $K_+=0.382$,
Theorem~\ref{thm:threshold} gives $\gc\approx2.04$, so no point blows up for
$\gamma\ge2.04$. Between these two values, the outcome depends on the whole profile and
on $c$, and we compute it numerically.

\begin{figure}[t]
	\centering
	\includegraphics[width=0.49\textwidth]{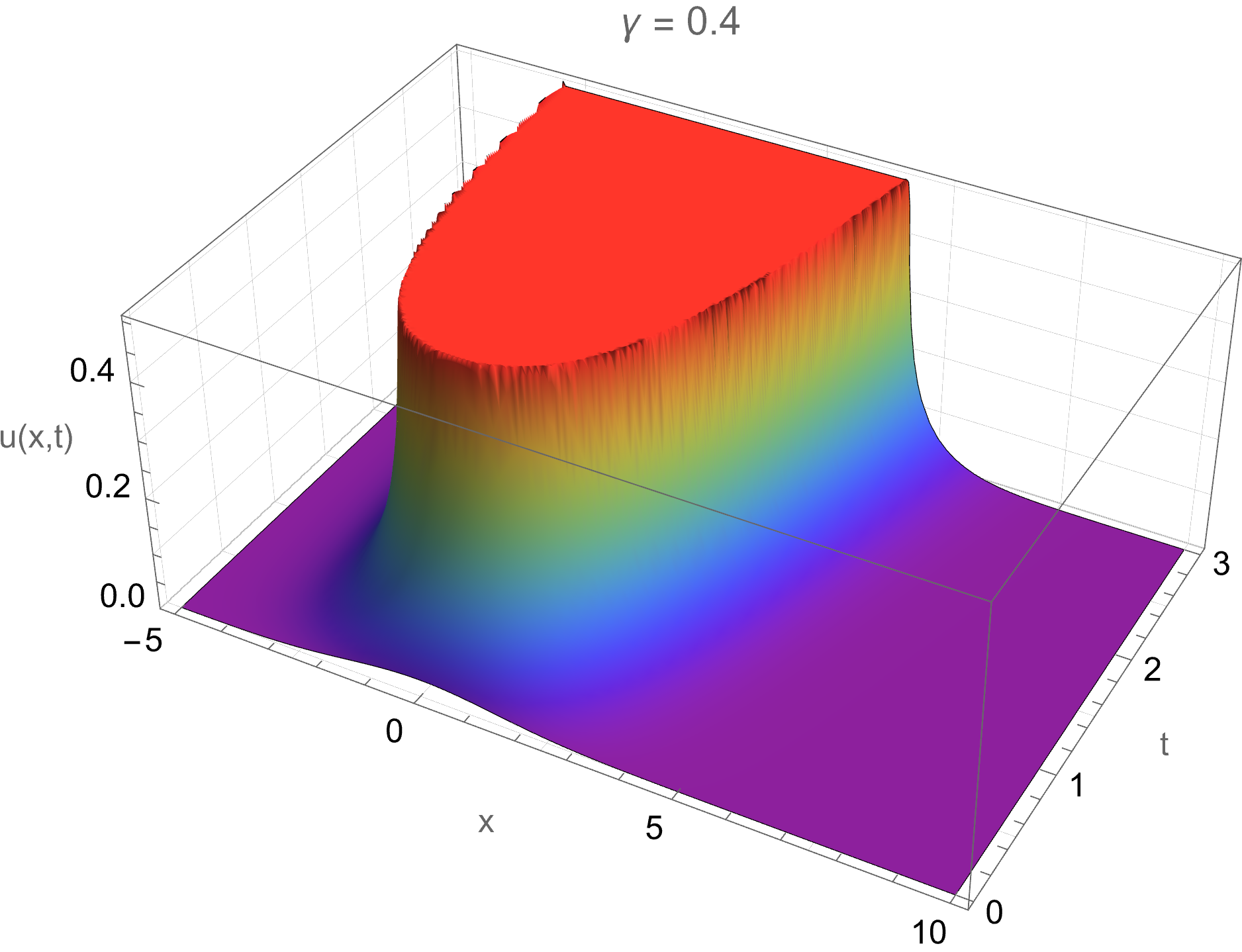}\hfill
	\includegraphics[width=0.49\textwidth]{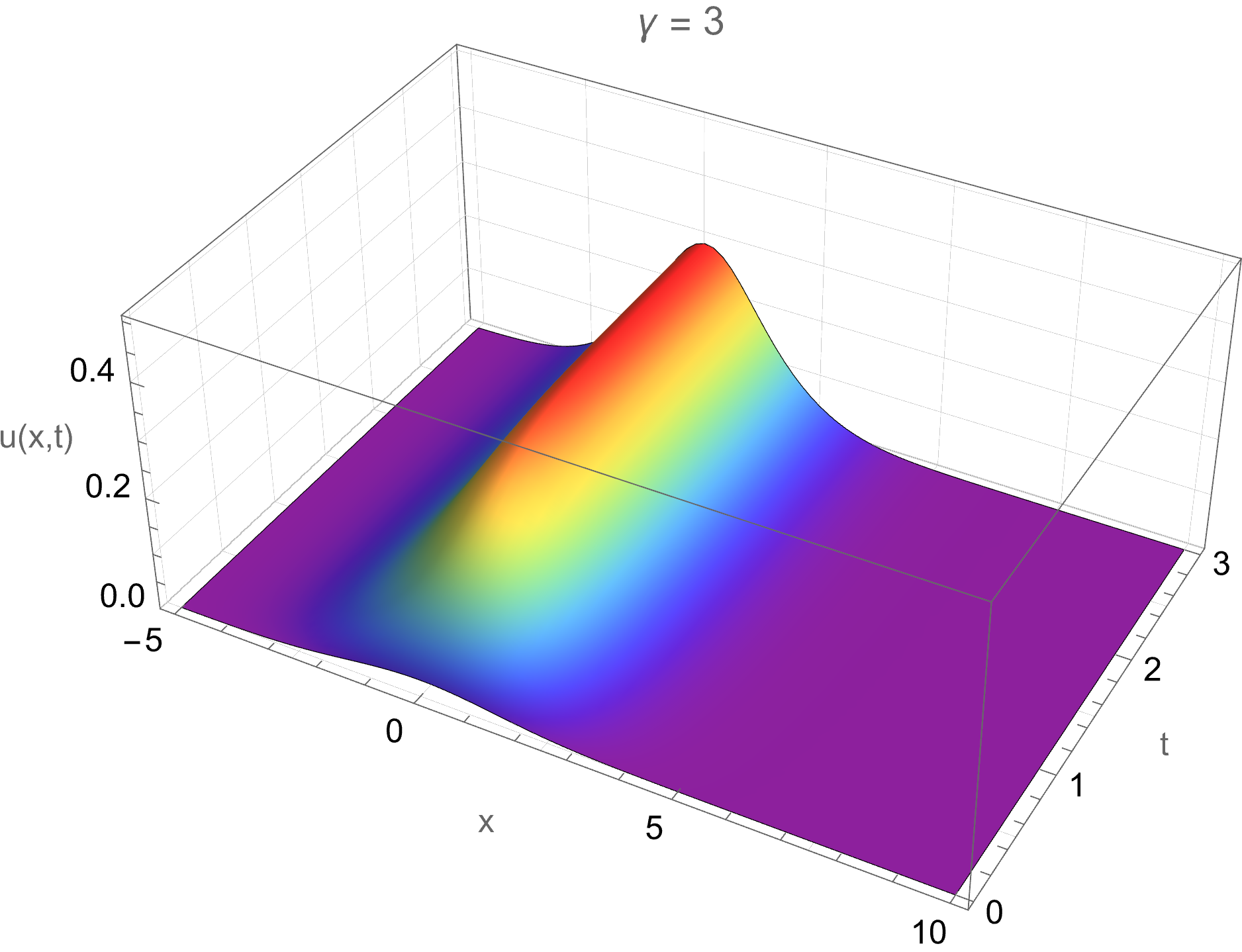}
\caption{Solution $u(x,t)$ of~\eqref{eq:model} for the Gaussian initial
	data~\eqref{eq:gauss}, with $\gamma=0.4$ (left) and $\gamma=3$ (right). In the left panel,
	values are cut off at a fixed level, so the flat region marks the points where the solution blows up.}
	\label{fig:surfaces}
\end{figure}

Figure~\ref{fig:surfaces} shows that for $\gamma=0.4$, blow-up starts at $x\approx0.1$,
$t^*\approx0.55$. By $t=3$, blow-up has occurred between $x\approx-2.8$ and $x\approx4.1$,
and over all time it occurs between $x\approx-3.2$ and $x\approx5.9$. The region extends farther downstream than upstream, because $K$ is
transported at speed $c$. For $\gamma=3$, the solution remains bounded everywhere and
drifts downstream. The set of points that blow up depends on the
initial profile and on $c$, not only on the local initial values.

\section{Applications}\label{sec:application}

We present three applications of the results of Section~\ref{sec:riccati}. In
Sections~\ref{sec:shock} and~\ref{sec:gel}, shock formation and gelation lead to the equation~\eqref{eq:reducedmodel} with constant $K$. In
Section~\ref{sec:friction}, the equation~\eqref{eq:model}, including transport,
follows from a heat balance.
\subsection{Shock formation under linear damping}\label{sec:shock}

As they propagate, nonlinear waves steepen while friction damps them. A finite
compression therefore forms a shock only if it steepens faster than friction can
smooth it out.
For the damped Burgers equation $\varphi_t+\varphi\varphi_x=-\nu\varphi$, with
$\varphi(x,0)=\varphi_0(x)$, $\varphi_0\in C^1$ and $\varphi_0'$ bounded below, the slope
$q=\varphi_x$ gives $\frac{dq}{dt}=-q^2-\nu q$ along each characteristic
$\frac{dx}{dt}=\varphi$. With $q=-e^{-\nu t}v$, we obtain
\[
\frac{dv}{dt}=e^{-\nu t}v^2,
\]
which has the form of the reduced equation~\eqref{eq:reducedmodel} with $\lambda=2$,
$K=0$ and $\gamma=\nu$, so the results of Section~\ref{sec:riccati} apply along each
characteristic, with initial value $A=-\varphi_0'(x_0)$, the initial compression.
Proposition~\ref{prop:power} with $n=1$ gives $\gc=A$. Therefore, a shock forms if and
only if the largest initial compression $A_{\max}=-\min\varphi_0'$ exceeds the damping
$\nu$, and it appears at the time \[t^*=-\frac{1}{\nu}\ln\bigl(1-\frac{\nu}{A_{\max}}\bigr)\]
given by~\eqref{eq:threshold}. The breaking of Burgers waves along characteristics is classical~\cite{Wh74,La64}, and the
criterion with linear damping appears, for example, in~\cite{Ch22}.

\subsection{Gelation with decaying reactivity}\label{sec:gel}
In polymerization and particle aggregation, clusters merge into ever larger clusters,
and the system may form a gel, a cluster of macroscopic size, in finite time. The
practical question is whether the system gels before its reactivity dies out.

In Smoluchowski's coagulation equation with the product kernel $r(t)\,ij$ and a monomer
input of rate $J(t)$, the second moment $M_2$ of the size distribution satisfies, before
gelation,
\[
\frac{dM_2}{dt}=r(t)M_2^2+J(t).
\]
For constant unit rates and $M_2(0)=0$, $M_2=\tan t$ and the system gels at
$t_g=\frac{\pi}{2}$~\cite{KM24,Le03}.
If the reactivity and the input decay together, $r=r_0e^{-\gamma t}$ and $J=J_0e^{-\gamma t}$,
for example through cooling or catalyst deactivation, this has the form of the reduced equation ~\eqref{eq:reducedmodel} with $\frac{\lambda}{2}=r_0$, $K=J_0$ and $A=M_2(0)$. With
input ($J_0>0$), Theorem~\ref{thm:threshold} shows that the system gels if and only if
$\gamma<\gc$. Without input ($J_0=0$), Proposition~\ref{prop:power} gives
$\gc=r_0M_2(0)$, and as $\gamma\to0$ the gelation time tends to the classical value
$\frac{1}{r_0M_2(0)}$.

\subsection{Frictional heating under a moving band}\label{sec:friction}
Frictional contacts are a classical source of self-excited
behaviour~\cite{Ib94a,Ib94b,PS90,Si10}. We now show that equation~\eqref{eq:model}
arises when classical models of frictional heating at a sliding contact are combined
with a viscoelastic model of the tension of an axially moving band~\cite{WM90}, under
the assumptions listed below.

\paragraph{Physical model.}
A band moves at constant speed $c$ over a fixed curved guide of radius $R$. The
coordinate $x$ is fixed in the laboratory, $\xi=x-ct$ labels the band material in
contact with the guide at $x$, and $u(x,t)$ is the temperature rise of the guide surface
at $x$. The model rests on five assumptions from the literature.
\begin{enumerate}[label=(\alph*),leftmargin=2.2em]
	\item \emph{Tension relaxation.} The band is modeled as a Maxwell element, as for
	conveyor belts~\cite{PC21}. Tensioned at $t=0$ and then held at constant length, its
	tension relaxes as $T(t)=T_0e^{-\gamma t}$, where $\gamma=1/\tau_r$ is the inverse
	relaxation time. Such exponential relaxation has been measured for polyester drive
	belts~\cite{Ka66}.
	\item \emph{Contact pressure.} A band under tension $T$ on a guide of radius $R$ presses
	on it with a normal force $T/R$ per unit length. We neglect variation in tension along the contact, which is small when the friction coefficient and wrap angle are small.
	\item \emph{Frictional heating.} At sliding speed $c$, friction produces a power
	$\mu\,(T/R)\,c$ per unit length, of which a fraction $\theta$ enters the
	guide~\cite{Ja42}.
	\item \emph{Temperature-dependent friction.} Nosko~\cite{No15} studies sliding
	instability with a friction coefficient that depends linearly on temperature. We keep
	the next term, $\mu=\mu_0(\xi)+\mu_1u+\mu_2u^2$ with $\mu_2>0$, which describes friction
	that rises increasingly fast with temperature, as measured for polycarbonate on
	aluminum above its glass transition~\cite{CLA23}. This law does not describe material pairs whose friction
	decreases with temperature~\cite{DC05}. A quadratic approximation
	of the source is also used in thermal explosion theory, where it leads to explicit
	solutions~\cite{BGW77}.
	\item \emph{Adiabatic guide.} The surface layer of the guide has heat capacity $m$ per
	unit length. Conduction along $x$ and heat losses are neglected. This is the adiabatic
	approximation of thermal explosion theory~\cite{FK69}.
\end{enumerate}

\paragraph{Derivation.}
Energy conservation at the point $x$ of the guide gives
\begin{equation}
	m\,u_t=\theta\,c\,\frac{T_0}{R}\,e^{-\gamma t}\Bigl[\mu_0(x-ct)+\mu_1u+\mu_2u^2\Bigr].
	\label{eq:heat}
\end{equation}
Set $\mathcal{A}=\frac{\theta cT_0}{mR}$ and $w=u+\frac{\mu_1}{2\mu_2}$. Completing the square gives
\begin{equation}
	w_t=e^{-\gamma t}\Bigl[\frac{\lambda}{2}w^2+K(x-ct)\Bigr],\qquad
	\lambda=2\mathcal{A}\mu_2,\qquad K(\xi)=\mathcal{A}\Bigl(\mu_0(\xi)-\frac{\mu_1^2}{4\mu_2}\Bigr),
	\label{eq:heatreduced}
\end{equation}
which is the reduced equation~\eqref{eq:reducedmodel}. Multiplying by $e^{\gamma t}$ and
applying $\partial_t+c\,\partial_x$, which annihilates $K(x-ct)$, gives
$w_{tt}+c\,w_{xt}+\gamma w_t=\lambda w(w_t+cw_x)e^{-\gamma t}$, which is
equation~\eqref{eq:model} for $w$. With $u$ in kelvin and $t$ in seconds, $\mathcal{A}$
and $K$ are in K\,s$^{-1}$ and $\lambda$ in K$^{-1}$s$^{-1}$, so every term
of~\eqref{eq:heatreduced} is in K\,s$^{-1}$.

In this model, $c$ is the band speed, $\gamma$ the tension-relaxation rate and $\lambda$
the strength of the thermal feedback. The damping term $\gamma w_t$ comes from the
relaxing tension, not from viscosity. The initial rate $u_t(x,0)$ is not free: it is
fixed by~\eqref{eq:heat} at $t=0$.

\paragraph{Meaning of the threshold.}
For a uniform band, $\mu_0$ and therefore $K$ are constant. When $K>0$,
Theorem~\ref{thm:threshold} applies at every point $x$, with
$A=u(x,0)+\frac{\mu_1}{2\mu_2}$, and $\tau^*$ is the adiabatic induction time at constant
tension, as in thermal explosion theory~\cite{FK69}. If the guide starts at the ambient
temperature, $u(x,0)=0$ and $A=\frac{\mu_1}{2\mu_2}$. The condition $K>0$ reads $\mu_0>\frac{\mu_1^2}{4\mu_2}$, that is, the friction law
$\mu_0+\mu_1u+\mu_2u^2$ has no real root. The compatibility condition of
Section~\ref{sec:riccati} holds automatically in this model, since $K$ comes from the
friction law. The condition $\gamma<\gc$ reads
$\tau_r>\tau^*$: the model predicts thermal runaway if and only if the tension relaxes
more slowly than runaway develops at full tension. For a purely elastic band
($\gamma\to0$), it predicts runaway whenever $K>0$. If $\mu_0$ varies along the band, for
example through a patch of higher friction, $K$ varies and Section~\ref{sec:surfaces}
applies.

Blow-up is a prediction of the adiabatic model, not a statement that the temperature
becomes infinite. It marks the onset of thermal runaway~\cite{BP07}. In a real guide,
conduction, heat losses, melting, wear, or a leveling-off of the friction curve bound
the temperature before the predicted blow-up time.

\section{Discussion: physical interpretation and limitations}\label{sec:discussion}
The threshold $\gc=\frac{1}{\tau^*}$ has a simple physical reading. $\tau^*$ is the blow-up time without damping ($\gamma=0$), and $\frac{1}{\gamma}$ is the total
rescaled time that the damping allows, however long the real time. For time-dependent
damping, $\tau_\infty$ plays this role (Theorem~\ref{thm:timedamping}). Blow-up occurs if and only if $\tau^*<\frac{1}{\gamma}$,
that is, if the nonlinear feedback is faster than the damping. In the physical model of
Section~\ref{sec:friction}, this is the race between thermal runaway and the relaxation
of the band's tension.

The results matter beyond the equation itself, for three reasons. First, the threshold
depends on the size of the disturbance, so linear stability analysis cannot predict it.
The explicit formula shows how large a disturbance must be to escape a given damping, and
gives a sharp boundary between finite and infinite growth. Second, the criterion
$\tau^*<\tau_\infty$ compares only two times, one set by the feedback and one by the
damping. It therefore applies to any system that reduces to the same form, as the shock
and gelation problems show, and it gives a design rule: in the friction model, runaway is
avoided if the tension relaxes faster than the induction time $\tau^*$. Third, the exact
blow-up times and the explicit bounds of Proposition~\ref{prop:bounds} provide benchmarks
for numerical schemes for blow-up problems, where errors are hard to measure.

The results have three limitations. First, the physical model of
Section~\ref{sec:friction} neglects heat conduction along the guide and heat losses,
which is the adiabatic approximation~\cite{FK69}. Adding a conduction or loss term breaks
the exact reduction, and the threshold is then expected to shift. Second, the friction
law $\mu=\mu_0+\mu_1u+\mu_2u^2$ is a truncated expansion of a measured
friction--temperature curve. Real curves level off or reverse at high
temperature~\cite{CLA23}, which bounds the runaway at a large but finite temperature,
as saturating friction laws bound self-excited motion~\cite{McM97,He12}. Third, the closed-form threshold requires $K$ to be the same on every characteristic.
When $K$ varies in space, the forcing at a fixed point changes with time.
Proposition~\ref{prop:bounds} still gives explicit bounds, and the solution must be
computed numerically only between them, as in Section~\ref{sec:surfaces}. These extensions, most
importantly heat conduction and heat losses in the guide, can be built on and compared
with the exact results of this paper.

\section{Conclusion}\label{sec:conclusion}

We have solved exactly, by the method of characteristics combined with a Riccati
reduction, a nonlinear damped transport-wave equation. For a natural class of initial data, the solution has a sharp, closed-form threshold, $\gc=\frac{1}{\tau^*}$, separating finite-time blow-up from global existence. This threshold depends on the finite size of the initial disturbance, so the
linearization of the equation cannot detect it. We extended it to time-dependent damping and to power-law feedback, gave explicit bounds
when $K$ varies in space, and illustrated the results numerically. We showed that the reduced equation describes exactly shock formation
under linear damping and gelation with decaying reactivity, and that the full equation
arises from a model of frictional heating under a moving viscoelastic band. Future work
should extend the analysis to saturating nonlinearities and to the heat conduction and
heat losses neglected in Section~\ref{sec:friction}.
\section*{Author contributions}
S.G. conceived the study, derived the theoretical results and wrote the manuscript.
S.B.K. carried out the numerical computations, prepared the figures, carried out the
literature search for the applications in Section~4 and contributed to revising the
manuscript. All authors read and approved the final version.
\bibliographystyle{unsrt}
\bibliography{refs}

\end{document}